\documentclass[11pt]{article}
\usepackage{latexsym} \usepackage{epsf}
\usepackage{a4}
\usepackage{amsfonts}
\renewcommand{\theequation}{\thesection.\arabic{equation}}
\newcounter{subequation}[equation]
\makeatletter

\expandafter\let\expandafter\reset@font\csname reset@font\endcsname

\def\subeqnarray{\arraycolsep1pt
    \def\@eqnnum\stepcounter##1{\stepcounter{subequation}%
        {\reset@font\rm(\theequation\alph{subequation})}}
\jot5mm     \eqnarray}

\makeatother

\def\tr{\mathop{\hbox{\rm tr}}\nolimits}

\def\be{\begin{equation}}
\def\ee{\end{equation}}
\def\bea{\begin{eqnarray}}
\def\eea{\end{eqnarray}}
\def\dd{\partial}
\def\half{\frac{1}{2}}

\def\one#1{#1^{\raise5pt\hbox{$\scriptstyle\!\!\!\!1$}}\,{}}
\def\two#1{#1^{\raise5pt\hbox{$\scriptstyle\!\!\!\!2$}}\,{}}

\def\II{\hbox{{1}\kern-.25em\hbox{l}}}
\makeatletter
\def\binrel@#1{\begingroup
  \setboxz@h{\thinmuskip0mu
    \medmuskip\m@ne mu\thickmuskip\@ne mu
    \setbox\tw@\hbox{$#1\m@th$}\kern-\wd\tw@
    ${}#1{}\m@th$}%
  \edef\@tempa{\endgroup\let\noexpand\binrel@@
    \ifdim\wdz@<\z@ \mathbin
    \else\ifdim\wdz@>\z@ \mathrel
    \else \relax\fi\fi}%
  \@tempa
}
\let\binrel@@\relax
\def\overset#1#2{\binrel@{#2}%
  \binrel@@{\mathop{\kern\z@#2}\limits^{#1}}}
\def\underset#1#2{\binrel@{#2}%
  \binrel@@{\mathop{\kern\z@#2}\limits_{#1}}}
\makeatother
\newfont{\bbd}{msbm10 scaled\magstep1}
\def\C{\hbox{\bbd C}}

\def\R{\hbox{\bbd R}}

\def\P{\hbox{\bbd P}}
\def\Q{{\mathbf Q}}

\def\RR{{\mathcal R}}

\def\RR{{\mathcal R}}

\newtheorem{prop}{Proposition}

\begin{document}

\vskip 1cm

\centerline{\LARGE\bf Factorization of the
$\mathrm{R}$-matrix }

\vskip 0.5cm 

\centerline{\LARGE\bf  and Baxter's Q-operator}

\vskip 1cm \centerline{\sc S.E. Derkachov} \vskip 1cm

\centerline{St.Petersburg Department of Steklov
Mathematical Institute of Russian Academy of Sciences,}
\centerline{Fontanka 27, 191023 St.Petersburg, Russia.}
\centerline{E-mail: {\tt derkach@euclid.pdmi.ras.ru}}

\vskip 2cm

{\bf Abstract.} The general rational solution of
the Yang-Baxter equation with the symmetry algebra
$s\ell(2)$ can be represented as the product of the 
simpler building blocks denoted as $\RR$-operators. The $\RR$-operators are
constructed explicitly and have simple structure. Using the 
$\RR$-operators we construct the two-parametric Baxter's 
Q-operator for the generic inhomogeneous XXX - spin chain.
In the case of homogeneous XXX-spin chain it is possible 
to reduce the general Q-operator to the much simpler 
one-parametric Q-operator.

\newpage

%%%%%%%%%%%%%%%%%%%%%%%%%%%%%%%%%%%%%%%%%%%%%%%%%%%%%%%%%%%%%%%%%%%%%%%%%%%%%%
{\small \tableofcontents}
\renewcommand{\refname}{References.}
\renewcommand{\thefootnote}{\arabic{footnote}}
\setcounter{footnote}{0}
\setcounter{equation}{0}

\renewcommand{\theequation}{\thesection.\arabic{equation}}
\setcounter{equation}{0}

\section{Introduction}
\setcounter{equation}{0}

The Yang-Baxter equation and its solutions play a key role
in the theory of the completely integrable quantum
models~\cite{KS,J,Dr,rep}. The general $s\ell(2)$-invariant 
solution of the Yang-Baxter equation (R-matrix) is the operator
$\mathbb{R}(u)$ acting in a tensor product of two 
$s\ell(2)$ lowest weights modules
$V_{\ell_1}\otimes V_{\ell_2}$. 
The Yang-Baxter equation is reduced to the simpler 
defining equation for the R-matrix~\cite{KRS}
$$
\R_{12}(u-v)\mathrm{L}_1(u)\mathrm{L}_2(v) =
\mathrm{L}_2(v)\mathrm{L}_1(u)\R_{12}(u-v)
$$
where $\mathrm{L}(u)$ is the Lax operator. 
We suggest the natural factorized expression for the
general R-matrix. It can be represented as the product of
the simple building blocks -- $\RR$-operators~\cite{Der}. The main
idea is very simple. The Lax operator depends on two
parameters: the spin of representation $\ell$ and the spectral
parameter $u$. It is useful to change to other parameters 
$u_{+} = u+\ell$ and $u_{-} = u-\ell$ and extract
the operator of permutation $\P_{12}$ from the R-matrix
$\R_{12}=\P_{12}\check{\R}_{12}$. The defining equation for
the operator $\check{\R}_{12}$ has the form
$$
\check{\R}_{12}\cdot
\mathrm{L}_1(u_{+},u_{-})\mathrm{L}_2(v_{+},v_{-}) =
\mathrm{L}_1(v_{+},v_{-})\mathrm{L}_2(u_{+},u_{-})\cdot
\check{\R}_{12}.
$$
The operator $\check{\R}_{12}$ interchanges simultaneously 
$u_{+}$ with $v_{+}$ and $u_{-}$ with $v_{-}$ in the product of two
Lax-operators. Let us perform this operation in two steps.
In the first step we interchange the parameters $u_{-}$ with
$v_{-}$ only. The parameters $u_{+}$ and $v_{+}$ remain the same.
In this way one obtains the natural defining equation for
the $\RR^{-}$-operator
$$
\RR^{-}_{12}\cdot 
\mathrm{L}_1(u_{+},u_{-})\mathrm{L}_2(v_{+},v_{-}) =
\mathrm{L}_1(u_{+},v_{-})\mathrm{L}_2(v_{+},u_{-})\cdot
\RR^{-}_{12}
\ ;\ \RR^{-}_{12} = \RR^{-}_{12}(u_{+},u_{-}|v_{-}).
$$
In the case when $u_{-}=v_{-}=v$ there is no interchange of 
parameters so that it is naturally to expect that 
the operator $\RR^{-}_{12}(u_{+},u_{-}|v_{-})$ 
is reduced to the unit operator $ \RR^{-}_{12}(u_{+},v|v) = \II $.
In the second step we interchange 
$u_{+}$ with $v_{+}$ but the parameters $u_{-}$ and $v_{-}$ 
remain the same.
The defining equation for the $\RR^{+}$-operator is
$$
\RR^{+}_{12}\cdot 
\mathrm{L}_1(u_{+},u_{-})\mathrm{L}_2(v_{+},v_{-}) =
\mathrm{L}_1(v_{+},u_{-})\mathrm{L}_2(u_{+},v_{-})\cdot
\RR^{+}_{12}
\ ;\ \RR^{+}_{12} = \RR^{+}_{12}(u_{+}|v_{+},v_{-}).
$$
In the case when $u_{+}=v_{+}=u$ there should be the 
similar degeneracy $\RR^{-}_{12}(u|u,v_{-}) = \II.$
These equations appear much simpler then the initial defining
equation for the R-operator and their solution can be
obtained in a closed form. Finally, we construct the
composite object - the R-matrix from the simplest building
blocks - the $\RR$-operators 
$$
\R_{12}(u_{+},u_{-}|v_{+},v_{-}) = 
\P_{12}\RR^{+}_{12}(u_{+}|v_{-},u_{-})
\RR^{-}_{12}(u_{+},u_{-}|v_{-}).
$$ 
There are two points of degeneracy $u_{-}=v_{-}=v$ and $u_{+}=v_{+}=u$ where 
the operator $\R_{12}$ is reduced to a single $\RR$-operator
$$
\R_{12}(u_{+},v|v_{+},v)= 
\P_{12}\RR^{+}_{12}(u_{+}|v_{+},v)
\ ;\ \R_{12}(u,u_{-}|u,v_{-})= 
\P_{12}\RR^{-}_{12}(u,u_{-}|v_{-})
$$
The detailed discussion of the R-matrix and its factorization is 
given in the Section 2.

The next natural step it to use the general 
operator $\R_{12}(u_{+},u_{-}|v_{+},v_{-})$
as building block in construction of Baxter's Q-operator.
In the case of the generic inhomogeneous periodic
XXX spin chain the transfer matrix $\mathbf{t}(u)$ is 
constructed as follows
$$
\mathbf{t}(u)=
\tr\mathrm{L}_1(u+\delta_1)\cdot\mathrm{L}_1(u+\delta_2)
\cdots\mathrm{L}_N(u+\delta_N)
$$ 
The most general transfer matrix $\mathbf{Q}(u,\ell_{0})$ 
is constructed in a similar manner from the 
operators $\R_{k0}$
$$
\mathbf{Q}(u,\ell_{0}) = 
\tr_{\mathrm{V}_{0}}\R_{10}(u+\delta_1)\R_{20}(u+\delta_2)
\cdots\R_{N0}(u+\delta_N)
$$
The operator $\mathbf{Q}(u,\ell_{0})$ depends on two parameters: 
the spectral parameter $u$ and the spin in the auxiliary space 
$\mathrm{V}_{0}=\mathrm{V}_{\ell_0}$.
It is useful to change to other parameters $u_1=1+u-\ell_0$
and $u_2 = u+\ell_{0}$ such that 
$\mathbf{Q}(u,\ell_{0})= \mathbf{Q}(u_1|u_2)$.
The operator $\mathbf{Q}(u_1|u_2)$ has the following properties:
\begin{itemize}
\item the operator  $\mathbf{Q}(u_1|u_2)$ is $s\ell(2)$-invariant
\item commutativity 
$$
\mathbf{Q}(u_1|u_2)\cdot\mathbf{Q}(v_1|v_2)=
\mathbf{Q}(v_1|v_2)\cdot\mathbf{Q}(u_1|u_2)
\  ; \  \mathbf{Q}(u_1|u_2)\cdot\mathbf{t}(v)=
\mathbf{t}(v)\cdot\mathbf{Q}(u_1|u_2)
$$
\item  the operator  $\mathbf{Q}(u_1|u_2)$ obeys the Baxter's equation
with respect to $u_2$ 
$$
\mathbf{Q}(u_1|u)\cdot\mathbf{t}(u) = 
\Delta_{+}(u)\mathbf{Q}(u_1|u+1)+
\Delta_{-}(u)\mathbf{Q}(u_1|u-1)
\ ;\ \Delta_{\pm}(u)=(u+\delta_1\pm\ell_1)\cdots
(u+\delta_N\pm\ell_N)
$$
\item the operator  $\mathbf{Q}(u_1|u_2)$ obeys the Baxter's equation
with respect to $u_1$
$$
\mathbf{t}(u)\cdot\mathbf{Q}(u|u_2) =
\frac{\Delta_{+}(u-1)\Delta_{-}(u)}{\Delta_{-}(u-1)}
\mathbf{Q}(u-1|u_2) + \Delta_{-}(u)\mathbf{Q}(u+1|u_2)
$$
\end{itemize}
These properties allow to consider the 
operator $\mathbf{Q}(u_1|u_2)$ as two-parametric 
Baxter's Q-operator~\cite{Bax,Skl2}.
The proof of all these properties of the operator 
$\mathbf{Q}(u_1|u_2)$ is given in the Section~3.

In the case of homogeneous spin chain:
$\delta_k=0$ and $\ell_k=\ell$ the points of degeneracy  
for all operators $\R_{k0}$ coincide so that it is possible 
to remove half of the $\RR$-operators.
We obtain the following reductions of the two-parametric Q-operator:
at the first point of degeneracy $u_1=1-\ell$
$$
\mathbf{Q}_{-}(u)=\mathbf{Q}(1-\ell|u)=
\tr_{\mathrm{V}_{0}}
\P_{10}\RR^{-}_{10}(u_{+},u_{-}|0)\cdot
\P_{20}\RR^{-}_{20}(u_{+},u_{-}|0)\cdots
\P_{N0}\RR^{-}_{N0}(u_{+},u_{-}|0)
$$
and at the second point of degeneracy $u_2=\ell$
$$
\mathbf{Q}_{+}(u)=\mathbf{Q}(u|\ell)=
\tr_{\mathrm{V}_{0}}
\P_{10}\RR^{+}_{10}(u_{+}|1,u_{-})\cdot
\P_{20}\RR^{+}_{20}(u_{+}|1,u_{-})\cdots
\P_{N0}\RR^{+}_{N0}(u_{+}|1,u_{-})
$$
As the direct consequence of the equations for the 
general two-parametric operator $\mathbf{Q}(u_1|u_2)$ 
we immediately derive the following properties of the 
operators $\mathbf{Q}_{+}(u)$ and $\mathbf{Q}_{-}(u)$
\begin{itemize}
\item operators $\mathbf{Q}_{\pm}(u)$ are $s\ell(2)$-invariant
\item commutativity
$$
\mathbf{Q}_{\pm}(u)\cdot\mathbf{Q}_{\pm}(v)
=\mathbf{Q}_{\pm}(v)\cdot\mathbf{Q}_{\pm}(u)
\ ;\ \mathbf{Q}_{+}(u)\cdot\mathbf{Q}_{-}(v)
=\mathbf{Q}_{-}(v)\cdot\mathbf{Q}_{+}(u)
\  ; \  \mathbf{Q}_{\pm}(u)\cdot\mathbf{t}(v)=
\mathbf{t}(v)\cdot\mathbf{Q}_{\pm}(u)
$$
\item Baxter equation for the $\mathbf{Q}_{-}(u)$
$$
\mathbf{Q}_{-}(u)\cdot\mathbf{t}(u) = 
\Delta_{+}(u)\mathbf{Q}_{-}(u+1)+
\Delta_{-}(u)\mathbf{Q}_{-}(u-1)
\ ;\ \Delta_{\pm}(u)=(u\pm\ell)^N
$$
\item Baxter equation for the $\mathbf{Q}_{+}(u)$
$$
\mathbf{t}(u)\cdot\mathbf{Q_{+}}(u) =
\frac{\Delta_{+}(u-1)\Delta_{-}(u)}{\Delta_{-}(u-1)}
\mathbf{Q}_{+}(u-1) + \Delta_{-}(u)\mathbf{Q}_{+}(u+1).
$$
\end{itemize}
There exists the natural generalization of the operators 
$\mathbf{Q}_{\pm}(u)$ to the case of generic 
inhomogeneous periodic XXX spin chain: we use the local 
operators $\RR_{k0}^{\pm}$ as building blocks 
for the operator $\mathbf{Q}_{\pm}(u)$.
These operators $\mathbf{Q}_{\pm}(u)$ obey 
the Baxter equations but they are not $s\ell(2)$-invariant 
and the commutation relations between 
$\mathbf{Q}_{\pm}(u)$ and between $\mathbf{Q}_{\pm}(u)$ 
and the transfer matrix $\mathbf{t}(u)$ are more complicated.
The two-parametric operator $\mathbf{Q}(u_1|u_2)$
can be factorized on the product of these operators
$$
\mathbf{Q}(u_1|u_2) = 
\mathbf{Q}_{+}(u_1)\cdot \P \cdot 
\mathbf{Q}_{-}(u_2)
$$
where $\P$ is the operator of cyclic shift.

The first explicit construction of the Q-operator 
using the general transfer matrix built from 
universal R-matrices was given by A.Yu.Volkov~\cite{Volkov} 
in context of some simple q-deformed model.
The very idea that the transfer matrix built from 
universal R-matrices can serve as a Baxter's Q-operator 
probably belongs to E.K.Sklyanin~\cite{Skl1,Skl2}. 
(The last sentence in fact coincides with footnote in Volkov's paper.)
Using universal R-matrix for the $U_q(\hat{s\ell}_2)$ affine algebra 
V.Bazhanov, S.Lukyanov and A.Zamolodchikov~\cite{BLZ} 
constructed Q-operator for quantum KdV model.
The general algebraic scheme how to derive the 
algebraic relations between different objects 
of QISM was formulated in the paper of A.Antonov and B.Feigin~\cite{AF}.   
The Baxter Q-operators was constructed 
for different models in the papers~\cite{PG,KSS,KP,RW,Z}.  
The factorization of the R-matrix and Baxter's 
Q-operator used in the present paper is 
very similar to the ones obtained in the context of the chiral 
Potts model~\cite{BS,Tar,BOU}.  

The presentation is organized as follows. In Section 2 we
collect the standard facts about the algebra~$s\ell(2)$ and
its representations. 
Next we consider the defining relation for the general
R-matrix, i.e. the solution of the Yang-Baxter equation
acting on tensor products of two arbitrary representations.
We introduce the natural defining equations for the
$\RR$-operators and show that the general R-matrix can be
represented as the product of these much simpler operators. 
In the Section 3 we construct Baxter's Q-operator for the 
generic inhomogeneous periodic XXX-spin chain.
In the Section 4 we consider the reduction to the case 
of the homogeneous spin chain. 
In the Section 5 we discuss the generalization of the 
operators $\mathbf{Q}_{\pm}(u)$ to the case of generic 
inhomogeneous spin chain and the factorization of the 
two-parametric operator $\mathbf{Q}(u_1|u_2)$.
Finally, in Section 6 we summarize.

\section{The general $s\ell(2)$-invariant R-matrix}
\setcounter{equation}{0}

The Lie algebra $s\ell(2)$ has three generators $\mathbf{S}\ ,\
\mathbf{S}_{\pm}$
$$
\left[\mathbf{S},\mathbf{S}_{\pm}\right] = \pm
\mathbf{S}_{\pm} \ ,\
\left[\mathbf{S}_{+},\mathbf{S}_{-}\right] = 2 \mathbf{S}\
$$
the central element (Casimir operator $\mathbf{C}_2$) being
$$
\mathbf{C}_2 = \mathbf{S}^2-\mathbf{S} +
\mathbf{S}_{+}\mathbf{S}_{-}.
$$
The Verma module $\mathbf{V}_{\ell}$ is the generic lowest
weight $s\ell(2)$-module with the lowest weight $\ell \in
\C$ and Casimir $\mathrm{C}_2 = \ell(\ell-1)$. As a linear
space $\mathbf{V}_{\ell}$ is spanned by the basis
$\left\{\mathbf{v}_k\right\}_{k=0}^{\infty}$
$$
\mathbf{v}_k = \mathbf{S}_{+}^k\mathbf{v}_{0} \ ,\
\mathbf{S}\mathbf{v}_k = (\ell+k)\mathbf{v}_k \ ,\
\mathbf{S}_{-}\mathbf{v}_k = -k(2\ell+k-1)\mathbf{v}_k
$$
where the vector $\mathbf{v}_0$ is the lowest weight
vector: $\mathbf{S}_{-}\mathbf{v}_0 = 0\ ,\
\mathbf{S}\mathbf{v}_0 = \ell \mathbf{v}_0$. The module
$\mathbf{V}_{\ell}$ is irreducible, except for $\ell = -
\frac{n}{2}$ where $ n \in \{0, 1, 2, 3 \cdots\}$, when
there exists an $(n+1)$-dimensional invariant subspace
$\mathbf{V}_{n}\subset\mathbf{V}_{\ell}$ spanned by
$\left\{\mathbf{v}_k\right\}_{k=0}^{n}$. We shall
rely extensively on the explicite representation $\mathrm{V}_{\ell}$ of
$s\ell(2)$ as the space $\C[z]$ of
polynomials in $z$ spanned by monomials 
$\left\{z^k\right\}_{k=0}^{\infty}$. Then the lowest weight
vector is polynomial $v_0 = 1$ and 
the action of $s\ell(2)$ in
$\mathrm{V}_{\ell}$ is given by the first-order
differential operators: 
\begin{equation}
\mathrm{S} = z\dd + \ell \ ,\
\mathrm{S}_{-} = -\dd\ ,\ \mathrm{S}_{+} = z^2\dd + 2\ell z.
\label{diff}
\end{equation} 
The generating function for the basis
vectors can be calculated in closed form 
$$
\mathrm{e}^{\lambda \mathrm{S}_{+}}\cdot 1 = (1-\lambda
z)^{-2\ell} = \sum_{k=0}^{\infty} \frac{\lambda^k}{k!}\cdot
\left(2\ell\right)_k z^{k} \ ;\ \left(2\ell\right)_k \equiv
\frac{\Gamma(2\ell+k)}{\Gamma(2\ell)}
$$
This expression clearly shows that for generic $\ell\neq
-\frac{n}{2}$ the module $\mathrm{V}_{\ell}$ is an
irreducible $s\ell(2)$-module isomorphic to
$\mathbf{V}_{\ell}$, the isomorphism being given by
$\mathbf{v}_k \leftrightarrow \left(2\ell\right)_k\cdot
z^k$. For $\ell = -\frac{n}{2}$ where $ n \in \{0, 1, 2, 3
\cdots\}$, we have the finite sum instead of infinite
series so that there exists an invariant subspace
$\mathrm{V}_{n}\subset\mathrm{V}_{\ell}$ spanned by
$\left\{z^k\right\}_{k=0}^{n}$ which is isomorphic to
$\mathbf{V}_{n}$. For $\ell = -\frac{1}{2}$ one obtains the
two-dimensional invariant subspace $\mathrm{V}_{1}\sim
\C^2$ and the matrices of operators $\mathrm{S}\ ,\
\mathrm{S}^{\pm}$ in the basis
$
\mathbf{e}_1 = \mathrm{S}_{+}\cdot 1 = -z ,\ \mathbf{e}_2 =
1
$ have the standard form of generators
$\mathbf{s} , \mathbf{s}_{\pm}$ in the fundamental
representation
\begin{equation}
\mathbf{s}_{+} = \left(\begin{array}{cc}
0 & 1 \\
0 & 0\end{array}\right)\ \ ,\ \mathbf{s}_{-} =
\left(\begin{array}{cc}
0 & 0 \\
1 & 0\end{array}\right)\ \ ,\ \mathbf{s} =
\half\cdot\left(\begin{array}{cc}
1 & 0 \\
0 & -1\end{array}\right)\ \label{fun}
\end{equation}
Let $\mathrm{V}_{\ell_1}$, $\mathrm{V}_{\ell_2}$ and
$\mathrm{V}_{\ell_3}$ be lowest weight $s\ell(2)$-modules
and consider three operators $\R_{\ell_i\ell_j}(u)$
which are acting in $V_{\ell_i}\otimes V_{\ell_j}$. The
Yang-Baxter equation is the following three term relation
\be \R_{\ell_1\ell_2}(u-v)\R_{\ell_1\ell_3}(u)
\R_{\ell_2\ell_3}(v)=
\R_{\ell_2\ell_3}(v)\R_{\ell_1\ell_3}(u)
\R_{\ell_1\ell_2}(u-v) \label{YB} \ee We look for the general
$s\ell(2)$-invariant solution $\R_{\ell_1\ell_2}(u)$ of
this equation. 
The restriction of the operator $\R_{\ell , -\half}(u)$ to the space
$\mathrm{V}_{\ell}\otimes\C^2$ coincides up to
normalization and a shift of the spectral parameter with
the fundamental Lax-operator~\cite{rep,KRS}
$$
\mathrm{L}(u): \mathrm{V}_{\ell}\otimes\C^2 \to
\mathrm{V}_{\ell}\otimes\C^2
$$
It is (up to an additive constant) the Casimir operator 
$\mathbf{C}_2$ for the tensor product of representations
$\mathrm{V}_{\ell}\otimes\C^2$\cite{KRS}
$$
\mathrm{L}(u) \equiv u  + 2\cdot\mathrm{S}\otimes
\mathbf{s}+\mathrm{S}_{-}\otimes
\mathbf{s}_{+}+\mathrm{S}_{+}\otimes\mathbf{s}_{-} =
\left(\begin{array}{cc}
u+\mathrm{S} & \mathrm{S}_{-} \\
\mathrm{S}_{+} & u-\mathrm{S}
   \end{array}\right) = \left(\begin{array}{cc}
u+\ell+z\dd & -\dd \\
z^2\dd + 2\ell z& u-\ell-z\dd
   \end{array}\right)
$$
where $\mathbf{s}, \mathbf{s}^{\pm}$ are the generators in
the fundamental representation~(\ref{fun}) and $\mathrm{S},
\mathrm{S}_{\pm}$ are the generators~(\ref{diff}) in the
generic representation~$\mathrm{V}_{\ell}$. The Lax
operator acts in the space $\C[z]\otimes\C^2$ and despite
of the compact notation $\mathrm{L}(u)$ depends actually on two
parameters: the spin $\ell$ and the spectral parameter $u$. We
shall use extensively the parametrization $u_{+}\equiv
u+\ell, u_{-}\equiv u-\ell$ and show all parameters
explicitly. There exists a very useful factorized
representation for the $\mathrm{L}$-operator~\cite{Skl1}
\begin{equation}
\mathrm{L}(u_{+},u_{-}) \equiv \left(\begin{array}{cc}
u_{+} +z\dd & -\dd \\
z^2\dd + (u_{+} -u_{-}) z& u_{-} -z\dd
   \end{array}\right) = \left(\begin{array}{cc}
1 & 0 \\
z& 1\end{array}\right)\ \left(\begin{array}{cc}
u_{+} -1 & -\dd \\
0& u_{-}\end{array}\right)\ \left(\begin{array}{cc}
1 & 0 \\
-z& 1\end{array}\right). \label{factor}
\end{equation}
We put $\ell_3 = -\frac{1}{2}$ in~(\ref{YB}) and
consider the restriction on the invariant subspace
$\mathrm{V}_{\ell_1}\otimes\mathrm{V}_{\ell_2}\otimes\C^2$.
In this way one obtains the defining equation for the
operator $\R_{12}(u)$~\cite{KRS,rep}
$$
\R_{12}(u-v)\mathrm{L}_1(u) \mathrm{L}_2(v) =
\mathrm{L}_2(v)\mathrm{L}_1(u) \R_{12}(u-v).
$$
The operator $\mathrm{L}_k$ acts nontrivially on the tensor
product $\mathrm{V}_{\ell_k}\otimes\C^2$ which is
isomorphic to $\C[z_k]\otimes\C^2$ and the operator
$\R_{12}(u)$ acts nontrivially on the tensor
product $\mathrm{V}_{\ell_1}\otimes\mathrm{V}_{\ell_2}$
which is isomorphic to $\C[z_1]\otimes\C[z_2] =
\C[z_1,z_2]$.
It is useful to extract the
operator of permutation $\P_{12}$
$$
\P_{12}\Psi(z_1,z_2) = \Psi(z_2,z_1)\ ;\
\Psi(z_1,z_2)\in\C[z_1,z_2]
$$
from the $\R$-operator $\R_{12}(u) =
\P_{12}\check{\R}_{12}(u)$ and solve the defining
equation for the $\check{\R}$-operator
$$
\check{\R}_{12}(u_{+},u_-|v_{+},v_-)
\mathrm{L}_{1}(u_+,u_-)\mathrm{L}_{2}(v_+,v_-)=
\mathrm{L}_{1}(v_+,v_-)\mathrm{L}_{2}(u_+,u_-)
\check{\R}_{12}(u_{+},u_-|v_{+},v_-)
$$
where $u_+ = u+\ell_1\ ,\ u_- = u-\ell_1\ ,\ v_+ =
v+\ell_2\ ,\ v_- = v - \ell_1$. The $\check{\R}$-operator
can be factorized into the product of the simpler elementary
building blocks - $\RR$-operators~\cite{Der}.
\begin{prop}
There exists operator $\RR^{+}_{12}$ which is the solution of
defining equations
\begin{equation} \RR^{+}_{12}
\mathrm{L}_{1}(u_+,u_-)\mathrm{L}_{2}(v_+,v_-)=
\mathrm{L}_{1}(v_+,u_-)\mathrm{L}_{2}(u_+,v_-)\RR^{+}_{12}
\label{F1}
\end{equation}
$$
\RR^{+}_{12} = \RR^{+}_{12}(u_+|v_+,v_-)\ ;\ \RR^{+}_{12}(u_+|v_+,v_-) =
\RR^{+}_{12}(u_++\lambda|v_+ +\lambda,v_- +\lambda).
$$
The system of equations~(\ref{F1}) for the operator $\RR^{+}_{12}$
is equivalent to the simpler system \be \RR^{+}_{12}\cdot
\left[\mathrm{L}_1(u_+,u_-)+ \mathrm{L}_2(v_+,v_-)\right] =
\left[ \mathrm{L}_1(v_+,u_-)+
\mathrm{L}_2(u_+,v_-)\right]\cdot \RR^{+}_{12} \ ;\ \RR^{+}_{12}\cdot z_1
= z_1 \cdot \RR^{+}_{12}.\label{F1def} \ee
These requirements fix the operator $\RR^{+}$ up to
an overall normalization constant. 
Fixing the normalization in a such way that 
$\RR^{+}: 1\mapsto 1$ we obtain
\begin{equation} 
\RR^{+}_{12}(u_{+}|v_+,v_-) = \frac{\Gamma(v_+-v_-)}{\Gamma(u_+-v_-)}
\frac{\Gamma(z_{21}\dd_2+u_+-v_-)}{\Gamma(z_{21}\dd_2+v_+-v_-)}\
;\ z_{21} = z_2-z_1.
\label{R+}
\end{equation}
\end{prop}
\begin{prop}
There exists operator $\RR^{-}_{12}$ which is the solution of
defining equations
\begin{equation}
\RR^{-}_{12} \mathrm{L}_{1}(u_+,u_-)\mathrm{L}_{2}(v_+,v_-)=
\mathrm{L}_{1}(u_+,v_-)\mathrm{L}_{2}(v_+,u_-) \RR^{-}_{12}
\label{F2}
\end{equation}
$$
\RR^{-}_{12} = \RR^{-}_{12}(u_+,u_-|v_-)\ ;\ \RR^{-}_{12}(u_+,u_-|v_-) =
\RR^{-}_{12}(u_+ +\lambda,u_- +\lambda|v_- +\lambda).
$$
The system of equations~(\ref{F2}) for the operator $\RR^{-}_{12}$
is equivalent to the simpler system \be \RR^{-}_{12}\cdot \left[
\mathrm{L}_1(u_+,u_-)+ \mathrm{L}_2(v_+,v_-)\right] =
\left[ \mathrm{L}_1(u_+,v_-)+
\mathrm{L}_2(v_+,u_-)\right]\cdot \RR^{-}_{12}\ ;\ \RR^{-}_{12}\cdot z_2
= z_2 \cdot \RR^{-}_{12}. \label{F2def} \ee
These requirements fix the operator $\RR^{-}_{12}$ up to
an overall normalization constant.
Fixing the normalization in a such way that 
$\RR^{-}: 1\mapsto 1$ we obtain
\begin{equation} 
\RR^{-}_{12}(u_+,u_-|v_-) = \frac{\Gamma(u_+-u_-)}{\Gamma(u_+-v_-)}
\frac{\Gamma(z_{12}\dd_1+u_+-v_-)}{\Gamma(z_{12}\dd_1+u_+-u_-)}\
;\ z_{12} = z_1-z_2.
\label{R-}
\end{equation}
\end{prop}
\begin{prop}
The operator $\check{\R}$ can be factorized in the following way
\begin{equation}
\check{\R}_{12}(u_{+},u_-|v_{+},v_-)=
\RR^{+}_{12}(u_+|v_+,u_-)\RR^{-}_{12}(u_+,u_-|v_{-}).
\label{Rfact}
\end{equation}
\end{prop}
Note that the relations~(\ref{F1def}),~(\ref{F2def}) 
are simply the rules of commutation of the $\RR$-operators with
$s\ell(2)$-generators written in a compact form. 
The $\RR$-operators change the spins of
$s\ell(2)$-representations 
$$ 
\RR^{+}_{12}(u_+|v_+,v_-) : V_{\ell_1}\otimes
V_{\ell_2} \to V_{\ell_1-\xi_{+}}\otimes V_{\ell_2+\xi_{+}} \
;\ \xi_{+} = \frac{u_{+}-v_{+}}{2} 
$$ 
$$
\RR^{-}_{12}(u_+,u_-|v_{-}) :
V_{\ell_1}\otimes V_{\ell_2} \to V_{\ell_1+\xi_2}\otimes
V_{\ell_2-\xi_{-}} \ ;\ \xi_{-} = \frac{u_{-}-v_{-}}{2}
$$ 
while the general R-matrix 
\begin{equation}
\R_{12}(u-v) =
\P_{12}\RR^{+}_{12}(u_+|v_+,u_-)\RR^{-}_{12}(u_+,u_-|v_{-})
\label{R}
\end{equation} 
appears automatically $s\ell(2)$-invariant
$
\left[\R_{12}(u) , \vec\mathbf{S}_1+\vec\mathbf{S}_2\right]= 0 
$
where $\vec\mathbf{S}_k\equiv(\mathbf{S}^{+}_k , 
\mathbf{S}_k , \mathbf{S}^{-}_k)$.

\section{Construction of the $\Q$-operator for
the generic inhomogeneous periodic XXX spin chain.}

\setcounter{equation}{0}

The transfer matrix $\mathbf{t}(u)$ 
for the generic inhomogeneous periodic
XXX spin chain is constructed as follows  
\begin{equation}
\mathbf{t}(u) = \tr
\mathrm{L}_{1}(u_1^{+},u_1^{-})\cdots
\mathrm{L}_{N}(u_N^{+},u_N^{-})
\ ;\ \mathrm{L}_{k}(u_k^{+},u_k^{-}) \equiv \left(\begin{array}{cc}
u_k^{+}+z_k\dd_{k} & -\dd_{k} \\
z_{k}^2\dd_{k} + (u_k^{+}-u_k^{-}) z_{k}& u_k^{-}-z_{k}\dd_{k}
   \end{array}\right),
\label{t}
\end{equation}
where $u_k^{\pm} = u+\delta_k\pm\ell_k$ and the trace is taken 
in the auxiliary space $\C^2$.
The most general transfer matrix $\mathbf{Q}(u_1|u_2)$ 
is constructed in a similar manner from the 
operators $\R_{k0}$
$$
\mathbf{Q}(u_1|u_2) = 
\tr_{\mathrm{V}_{0}}\R_{10}(u+\delta_1)\R_{20}(u+\delta_2)
\cdots\R_{N0}(u+\delta_N)
$$
where we use the parameters $u_1=1+u-\ell_0$
and $u_2 = u+\ell_{0}$ instead of the spectral 
parameter $u$ and the spin parameter $\ell_0$ 
in the auxiliary space $\mathrm{V}_{0}=\mathrm{V}_{\ell_0}$.
The explicit expression for the operator $\R_{k0}$ is~(\ref{R})
$$
\R_{k0}(u+\delta_k)=\frac{\Gamma(\ell_k+\ell_0-u-\delta_k)}
{\Gamma(\ell_k+\ell_0+u+\delta_k)}
\cdot\P_{k0}\cdot
\frac{\Gamma(z_{0k}\dd_0+2\ell_k)}
{\Gamma(z_{0k}\dd_0+\ell_k+\ell_0-u-\delta_k)} 
\frac{\Gamma(z_{k0}\dd_k+\ell_k+\ell_0+u+\delta_k)}
{\Gamma(z_{k0}\dd_k+2\ell_k)}
$$
and it is natural to simplify the notations: 
we omit the local parameters $\ell_k,\delta_k$ in the chain
and show the global parameters $u_1=1+u-\ell_0$
and $u_2 = u+\ell_{0}$ only 
\begin{equation}
\R_{k0}(u_1|u_2) =
\frac{\Gamma(\ell_k+1-u_1-\delta_k)}
{\Gamma(\ell_k+u_2+\delta_k)}
\cdot\P_{k0}\cdot
\frac{\Gamma(z_{0k}\dd_0+2\ell_k)}
{\Gamma(z_{0k}\dd_0+\ell_k+1-u_1-\delta_k)} 
\frac{\Gamma(z_{k0}\dd_k+\ell_k+u_2+\delta_k)}
{\Gamma(z_{k0}\dd_k+2\ell_k)}
\label{Rglob}
\end{equation}
The basic properties of the operator $\mathbf{Q}(u_1|u_2)$ 
are enumerated in the Introduction. 

The $s\ell(2)$-invariance: 
$
\left[\mathbf{Q}(u_1|u_2) ,\vec\mathbf{S}_1+\cdots+\vec\mathbf{S}_N\right]= 0
$
follows immediately from the $s\ell(2)$-invariance of operators 
$\R_{k0}$: 
$
\left[\R_{k0}(u) , \vec\mathbf{S}_k+\vec\mathbf{S}_0\right]=0 
$
and the cyclicity property of the trace.

The commutativity 
$
\left[\mathbf{Q}(u_1|u_2), \mathbf{Q}(v_1|v_2)\right]= 0
$
follows from the Yang-Baxter equation for the general R-matrix
$$
\R_{00^{\prime}}(u-v)\R_{k0}(u)\R_{k0^{\prime}}(v) = 
\R_{k0^{\prime}}(u)\R_{k0}(u)\R_{00^{\prime}}(u-v)
$$
where $\mathrm{V}_0=\mathrm{V}_{\ell_0}$ and 
$\mathrm{V}_{0^{\prime}}=\mathrm{V}_{\ell_{0^{\prime}}}$ 
are two auxiliary spaces and $\mathrm{V}_k=\mathrm{V}_{\ell_k}$ 
is the k-th quantum space.

The commutativity
$
\left[\mathbf{Q}(u_1|u_2), \mathbf{t}(v)\right] = 0
$
follows from the special 
case of the general Yang-Baxter relation
$$
\R_{k0}(u-v)\mathrm{L}_k(u)\mathrm{L}_0(v) =
\mathrm{L}_0(v)\mathrm{L}_k(u)\R_{k0}(u-v).
$$
All these formulae are standard and well known. 
The really nontrivial is the derivation of the Baxter equations. 
It is the consequence of the important properties of the $\RR$-operators -
triangularity relations.
\begin{prop}
The following triangularity relations hold for the 
operators $\RR^{-}_{12}\left(u_{+},u_{-}|0\right)$ and
$\RR^{+}_{12}\left(u_{+}|1,u_{-}\right)$ 
\begin{equation}
\mathbf{M}^{-1}_1\cdot\RR^{-}_{12}\left(u_{+},u_{-}|0\right)\cdot
\mathrm{L}_1\left(u_{+},u_{-}\right)
\cdot\mathbf{M}_2
= \label{trianglR-}
\end{equation}
$$
=\left(\begin{array}{cc} u_{+}\cdot\RR^{-}_{12}\left(u_{+}+1,u_{-}+1|0\right) &
-\RR^{-}_{12}\left(u_{+},u_{-}|0\right)\dd_1\\
0 & u_{-}\cdot\RR^{-}_{12}\left(u_{+}-1,u_{-}-1|0\right)
\end{array}\right),\
$$
\begin{equation}
\mathbf{M}^{-1}_1\cdot \mathrm{L}_2\left(u_{+},u_{-}\right)\cdot
\RR^{+}_{12}\left(u_{+}|1,u_{-}\right) 
\cdot\mathbf{M}_2
= \label{trianglR+}
\end{equation}
$$
= \left(\begin{array}{cc} 
\frac{u_-(u_+ -1)}{u_- -1}\cdot\RR^{+}_{12}\left(u_{+}-1|1,u_{-}-1\right) &
-\dd_1\RR^{+}_{12}\left(u_{+}|1,u_{-}\right)\\
0 & u_-\cdot\RR^{+}_{12}\left(u_{+}+1|1,u_{-}+1\right)
\end{array}\right),\
$$
where 
$$
\RR^{-}_{12}(u_+,u_-|0) = \frac{\Gamma(u_+-u_-)}{\Gamma(u_+)}
\frac{\Gamma(z_{12}\dd_1+u_+)}{\Gamma(z_{12}\dd_1+u_+-u_-)}
\ ;\ 
\RR^{+}_{12}(u_{+}|1,u_-) = \frac{\Gamma(1-u_-)}{\Gamma(u_+-u_-)}
\frac{\Gamma(z_{21}\dd_2+u_+-u_-)}{\Gamma(z_{21}\dd_2+1-u_-)},
$$
$$
\mathbf{M}_k \equiv \left(\begin{array}{cc}
1 & 0 \\
z_k& 1\end{array}\right).
$$
\end{prop}
We prove the triangularity relation for the operator $\RR^{-}_{12}$ and 
the proof for the operator $\RR^{+}_{12}$ is very similar.
We start directly from the defining equation~(\ref{F2}).
Using factorization~(\ref{factor}) of the Lax operator and
commutativity of $\RR^{-}_{12}$ and $z_2$ the defining equation for the
$\RR^{-}_{12}$-operator can be represented in the form
$$
\RR^{-}_{12} \left(\begin{array}{cc}
u_++z_1\dd_{1} & -\dd_{1} \\
z_{1}^2\dd_{1} + (u_+-u_-) z_{1}& u_--z_{1}\dd_{1}
   \end{array}\right)\left(\begin{array}{cc}
1& 0 \\
z_{2}& 1
   \end{array}\right)
\left(\begin{array}{cc}
1& -\dd_{2} \\
0& v_-
   \end{array}\right)=
$$
$$
= \left(\begin{array}{cc}
1 & 0 \\
z_1& 1\end{array}\right)\ \left(\begin{array}{cc}
u_+ -1 & -\dd_1 \\
0& v_-\end{array}\right)\ \left(\begin{array}{cc}
1 & 0 \\
-z_1& 1\end{array}\right) \left(\begin{array}{cc}
1& 0 \\
z_{2}& 1
   \end{array}\right)
\left(\begin{array}{cc}
1& -\dd_{2} \\
0& u_-
   \end{array}\right)\RR^{-}_{12}
$$
Next we transform all this in the following simple way
$$
\left(\begin{array}{cc}
1 & 0 \\
-z_1& 1\end{array}\right)\ \RR^{-}_{12} \left(\begin{array}{cc}
u_++z_1\dd_{1} & -\dd_{1} \\
z_{1}^2\dd_{1} + (u_+-u_-) z_{1}& u_--z_{1}\dd_{1}
   \end{array}\right)\left(\begin{array}{cc}
1& 0 \\
z_{2}& 1
   \end{array}\right)=
$$
$$
= v_-^{-1}\cdot \left(\begin{array}{cc}
u_+ -1 & -\dd_1 \\
0& v_-\end{array}\right)\ \left(\begin{array}{cc}
1 & 0 \\
-z_{12}& 1\end{array}\right) \left(\begin{array}{cc}
1& -\dd_{2} \\
0& u_-
   \end{array}\right) \RR^{-}_{12} \left(\begin{array}{cc}
v_-& \dd_{2} \\
0& 1
   \end{array}\right) =
$$
$$
=\left(\begin{array}{cc} (u_++z_{12}\dd_1)\cdot\RR^{-}_{12} & 
-\RR^{-}_{12}\dd_1\\
-v_-\cdot z_{12}\RR^{-}_{12} &
\frac{u_-(u_+-1)-v_-(u_--z_{12}\dd_1)}{u_+-v_--1-z_{12}\dd_1}\cdot\RR^{-}_{12}
\end{array}\right)\
$$
In the last transformation we use the explicit representation for the
$\RR^{-}_{12}$-operator~(\ref{R-}). 
The obtained matrix becomes triangular at the point
$v_{-}=0$ and using the explicit expression for 
the operator $\RR^{-}_{12}$ this matrix can be transformed 
to the form~(\ref{trianglR-}).
From the triangularity relations for the
$\RR$-operators immediately follow two triangularity
relations for the general operator $\R_{12}$.

\begin{prop}
The following triangularity relations for the 
operator $\R_{12}\left(u_{+},u_{-}|v_+,v_-\right)$ hold 
\begin{equation}
\mathbf{M}^{-1}_2\cdot
\R_{12}\left(u_{+},u_{-}|v_+ ,0\right)\cdot
\mathrm{L}_1\left(u_{+},u_{-}\right)
\cdot\mathbf{M}_2
= \label{trianglR1}
\end{equation}
$$
=\left(\begin{array}{cc} u_{+}\cdot
\RR_{12}\left(u_{+}+1,u_{-}+1|v_+ +1, 0\right) &
*** \\
0 & u_{-}\cdot\R_{12}\left(u_{+}-1,u_{-}-1|v_+ -1, 0\right)
\end{array}\right)\
$$
\begin{equation}
\mathbf{M}^{-1}_2\cdot
\mathrm{L}_1\left(u_{+},u_{-}\right)\cdot
\R_{12}\left(u_{+},u_-|1,v_{-}\right) 
\cdot\mathbf{M}_2
= \label{trianglR2}
\end{equation}
$$
= \left(\begin{array}{cc} 
\frac{u_-(u_+ -1)}{u_- -1}\cdot\R_{12}\left(u_{+}-1,u_- -1|1,v_{-}-1\right) &
***\\
0 & u_-\cdot\R_{12}\left(u_{+}+1,u_- +1|1,v_{-}+1\right)
\end{array}\right)\
$$
\end{prop}
The relation~(\ref{trianglR1}) is obtained from the 
relation~(\ref{trianglR-}) simply by multiplying 
with the operator $\P_{12}\RR^{+}_{12}(u_+|v_+, u_-)$ from 
the left and using the expression~(\ref{Rfact}) 
for the operator $\R_{12}$. 
The operator $\RR^+$ plays a passive role in 
this first relation.
The relation~(\ref{trianglR2}) is obtained from the 
relation~(\ref{trianglR+}) simply by multiplying 
with the operator $\P_{12}$ from the left and with the operator 
$\RR^{-}_{12}(u_+,u_-|v_-)$ from the right and 
using the expression~(\ref{Rfact}) for the operator $\R_{12}$. 
Now the operator $\RR^-$ plays a a passive role.

Let us go to the proof of the Baxter equation
$$
\mathbf{Q}(u_1|u)\cdot\mathbf{t}(u) = 
\Delta_{+}(u)\mathbf{Q}(u_1|u+1)+
\Delta_{-}(u)\mathbf{Q}(u_1|u-1)
\ ;\ \Delta_{\pm}(u)=(u+\delta_1\pm\ell_1)\cdots
(u+\delta_N\pm\ell_N)
$$
It is the direct consequence of the triangularity 
relation~(\ref{trianglR1}) and cyclicity of the trace.
Let us choose the first space in~(\ref{trianglR1})
as k-th quantum space and the second 
space as the auxiliary space.
We have in useful notations
$$
\mathbf{M}^{-1}_0\cdot\R_{k0}\left(u_1-v|u\right)\cdot
\mathrm{L}_k\left(u+\delta_k\right)
\cdot\mathbf{M}_0 = 
\left(\begin{array}{cc} u^{+}_k\cdot
\R_{k0}\left(u_1-v|u+1\right) &
*** \\
0 & u^{-}_k\cdot\R_{k0}\left(u_1-v|u-1\right)
\end{array}\right)\
$$
Multipying these equalities for $k=1, 2, \cdots N $ , 
taking the traces in auxiliary spaces $\C^2$ and 
$\mathrm{V}_0$ and using the cyclisity of the 
trace one obtains the equation
$$
\mathbf{Q}(u_1-v|u)\cdot\mathbf{t}(u) = 
\Delta_{+}(u)\mathbf{Q}(u_1-v|u+1)+
\Delta_{-}(u)\mathbf{Q}(u_1-v|u-1)
$$
The parameter $u_1$ is arbitrary so that 
we obtain the needed relation.
The Baxter equation with respect to parameter $u_2$ 
follows from the triangularity relation~(\ref{trianglR2})
and the derivation is very similar.

\section{Q-operator for the homogeneous periodic XXX spin chain}

\setcounter{equation}{0}

The operator $\R_{k0}(u_1|u_2)$~(\ref{Rglob}) has two 
points of degeneracy: $u_1 = 1-\delta_k-\ell_k$ and 
$u_2 = \ell_k-\delta_k$.
In the case of homogeneous spin chain:
$\delta_k=0$ and $\ell_k=\ell$, the degeneration points 
for all operators $\R_{k0}$ coincide so that it is possible 
to remove half of the $\RR$-operators in 
the two-parametric operator 
$$
\mathbf{Q}(u_1|u_2) = 
\tr_{\mathrm{V}_{0}}\R_{10}(u_1|u_2)\R_{20}(u_1|u_2)
\cdots\R_{N0}(u_1|u_2) 
$$
We obtain the following reductions of the 
two-parametric Q-operator:
at the first point of degeneracy $u_1=1-\ell$
$$
\mathbf{Q}_{-}(u)=\mathbf{Q}(1-\ell|u)=
\tr_{\mathrm{V}_{0}}
\P_{10}\RR^{-}_{10}(u_{+},u_{-}|0)\cdot
\P_{20}\RR^{-}_{20}(u_{+},u_{-}|0)\cdots
\P_{N0}\RR^{-}_{N0}(u_{+},u_{-}|0)=
$$
$$
=\frac{\Gamma^N(2\ell)}{\Gamma^N(\ell+u)}
\cdot\tr_{\mathrm{V}_{0}}\P_{10}
\frac{\Gamma\left(z_{10}\dd_1+u+\ell\right)}
{\Gamma\left(z_{10}\dd_1+2\ell\right)}\cdots
\P_{N0}
\frac{\Gamma\left(z_{N0}\dd_N+u+\ell\right)}
{\Gamma\left(z_{N0}\dd_N+2\ell\right)}
$$
and at the second point of degeneracy $u_2=\ell$
$$
\mathbf{Q}_{+}(u)=\mathbf{Q}(u|\ell)=
\tr_{\mathrm{V}_{0}}
\P_{10}\RR^{+}_{10}(u_{+}|1,u_{-})\cdot
\P_{20}\RR^{+}_{20}(u_{+}|1,u_{-})\cdots
\P_{N0}\RR^{+}_{N0}(u_{+}|1,u_{-})=
$$
$$
=\frac{\Gamma^N(1+\ell-u)}{\Gamma^N(2\ell)}\cdot
\tr_{\mathrm{V}_{0}}\P_{10}
\frac{\Gamma\left(z_{01}\dd_0+2\ell\right)}
{\Gamma\left(z_{01}\dd_0+1+\ell-u\right)}\cdots
\P_{N0}
\frac{\Gamma\left(z_{0N}\dd_0+2\ell\right)}
{\Gamma\left(z_{0N}\dd_0+1+\ell-u\right)}.
$$
As the direct consequence of the equations for the 
general two-parametric operator $\mathbf{Q}(u_1|u_2)$ 
we immediately obtain the corresponding properties of the 
operators $\mathbf{Q}_{+}(u)$ and $\mathbf{Q}_{-}(u)$ 
itemized in Introduction.
We construct the $\mathbf{Q}_{\pm}$-operator as the trace of the
products of $\RR^{\pm}$ operators in auxiliary space $V_0$. The whole
construction is pure algebraic. 
In this Section we shall derive the explicit formulae
for the action of the operator $\mathbf{Q}_{-}$ 
in the space of polynomials. The explicit expression for the 
second operator $\mathbf{Q}_{+}$ is more complicated 
and we shall not consider it here. 
We have the following expression for the operator $\mathbf{Q}_{-}$
$$
\mathbf{Q}_{-}(u) = 
\tr_{\mathrm{V}_{0}}\P_{10}\RR\left(z_{10}\dd_1\right)\cdots
\P_{N0}\RR\left(z_{N0}\dd_N\right)\ ;\ \RR(x)\equiv
\frac{\Gamma\left(x+u+\ell\right)}
{\Gamma\left(x+2\ell\right)}. 
$$
\begin{prop}
The action of the operator $\mathbf{Q}_{-}(u)$ on a polynomial 
$\Psi(z_1\cdots z_N)$ can be represented in the 
following equivalent forms
\begin{eqnarray}
\label{QL1}
&&\left[\mathbf{Q}_{-}(u)\Psi\right](z_1,\cdots z_N)=\\[2mm]
&&\ \ \ \ \ \ \left.
\RR(z_{01}\dd_0)\RR(z_{12}\dd_1)\RR(z_{23}\dd_2)\cdots
\RR(z_{N-1,N}\dd_N)\Psi(z_0,z_1,\cdots z_{N-1})
\right|_{z_0=z_N} \nonumber\\[4mm]
\label{QL2}
&&\left[\mathbf{Q}_{-}(u)\Psi\right](z_1,\cdots z_N)=\\[2mm]
&&\ \ \ \ \ \ \ \ \RR(t_1\dd_{t_1})\RR(t_2\dd_{t_2}),\cdots \RR(t_N\dd_{t_N})
\cdot \left.\Psi\left(t_1z_{N1}+z_1,t_2z_{12}+z_2\cdots
t_Nz_{N-1,N}+z_N\right)\right|_{t_k=1}\nonumber
\end{eqnarray}
\begin{equation}
\left[\mathbf{Q}_{-}(u)\Psi\right](z_1\cdots z_N) =
\frac{\Gamma^N(2\ell)}{\Gamma^N(\ell+u)\Gamma^N(\ell-u)}\cdot
\int^{1}_{0}\mathrm{d}\alpha_1 (1-\alpha_1)^{\ell-u-1}
\alpha_1^{\ell+u-1}\cdots 
\label{QL3}
\end{equation}
$$
\cdots 
\int^{1}_{0}\mathrm{d}\alpha_N 
(1-\alpha_N)^{\ell-u-1}\alpha_N^{\ell+u-1}
\Psi\left(\alpha_1z_{N1}+z_1 , \alpha_2z_{12}+z_2 \cdots
\alpha_Nz_{N-1,n}+z_N\right).
$$
The operator $\mathbf{Q}_{-}(u)$ maps polynomials 
in variables $z_1 \cdots z_N$ to
polynomials in variables $u,z_1 \cdots z_N$
$$
\mathbf{Q}_{-}(u) :\ \C[z_1\cdots z_N] \mapsto \C[u,z_1 \cdots z_N].
$$
\end{prop}
Let $z_0$ be the variable in the auxiliary space $V_0$
and let the operator $\mathbb{A}$ act in the tensor product
$\mathrm{V}_{0}\otimes\mathrm{V}_{1}\cdots\otimes \mathrm{V}_{N}$ 
and $\Psi(z_1\cdots z_N) \in
\mathrm{V}_{1}\cdots\otimes \mathrm{V}_{N}$. 
The trace of the operator $\mathbb{A}$ in auxiliary
space $\mathrm{V}_{0}=\C[z_0]$ can be calculated as follows
$$
\left[\left(\tr_{\mathrm{V}_{0}}
\mathbb{A}\right) \Psi\right](z_1\cdots z_N) =
\left.\sum_{m=0}^{+\infty} \frac{1}{m!} \dd_0^m \mathbb{A} \cdot z_0^m
\cdot\Psi(z_1\cdots z_N)\right|_{z_0=0}.
$$
In order to prove~(\ref{QL1}) it is useful to move all 
permutations to the right
$$ \P_{10}\RR(z_{10}\dd_1)\P_{20}\RR(z_{20}\dd_2)\cdots
\P_{N0}\RR(z_{N0}\dd_N) =
$$
$$
=\RR(z_{01}\dd_0)\RR(z_{12}\dd_1)\RR(z_{23}\dd_2)\cdots
\RR(z_{N-1,N}\dd_{N-1}) \cdot
\P_{10}\P_{20}\cdots\P_{N0}.
$$
Then we have
$$
\RR(z_{01}\dd_0)\RR(z_{12}\dd_1)\RR(z_{23}\dd_2)\cdots
\RR(z_{N-1,N}\dd_{N-1}) \cdot
\P_{10}\P_{20}\cdots\P_{N0}\cdot z_0^m
\cdot\Psi(z_1\cdots z_N) =  
$$
$$
=\RR(z_{01}\dd_0)\RR(z_{12}\dd_1)\RR(z_{23}\dd_2)\cdots
\RR(z_{N-1,N}\dd_{N-1})\cdot z_N^m \cdot\Psi(z_0,z_1\cdots z_{N-1})=
$$
$$
= z_N^m \RR(z_{01}\dd_0)\RR(z_{12}\dd_1)\RR(z_{23}\dd_2)\cdots
\RR(z_{N-1,N}\dd_{N-1})\Psi(z_0,z_1\cdots z_{N-1}).
$$
The result of the operation $\sum_{m=0}^{+\infty} \frac{1}{m!} \dd_0^m$ can
be calculated in closed form 
$$
\left[\mathbf{Q}_{-}(u)\Psi\right](z_1\cdots z_N) = 
\sum_{m=0}^{+\infty} \frac{1}{m!}
\dd_0^m z_N^m \left.\RR(z_{01}\dd_0)\RR(z_{12}\dd_1)\cdots
\RR(z_{N-1,N}\dd_{N-1})\Psi(z_0,z_1\cdots z_{N-1})\right|_{z_0=0}=
$$
$$
= \mathrm{e}^{z_N\dd_0}
\left.\RR(z_{01}\dd_0)\RR(z_{12}\dd_1)\RR(z_{23}\dd_2)\cdots
\RR(z_{N-1,N}\dd_{N-1})\Psi(z_0,z_1\cdots z_{N-1})\right|_{z_0=0}=
$$
$$
= \left.\RR(z_{01}\dd_0)\RR(z_{12}\dd_1)\RR(z_{23}\dd_2)\cdots
\RR(z_{N-1,N}\dd_{N-1})\Psi(z_0,z_1\cdots z_{N-1})\right|_{z_0=z_N}
$$
and we obtain the first representation~(\ref{QL1}).

The second formula is the simple consequence of the first one. Let us
represent the action of the operator $\RR(z_{k-1,k}\dd_{k-1})$ in the
form
$$
\RR(z_{k-1,k}\dd_{k-1})\Psi(z_{k-1}) = \left.\RR(t_k\dd_{t_k})
t_k^{z_{k-1,k}\dd_{k-1}} \Psi(z_{k-1})\right|_{t_k=1} =
$$
$$
=\left.\RR(t_k\dd_{t_k})\mathrm{e}^{-z_k\dd_{z_{k-1}}}
t_k^{z_{k-1}\dd_{k-1}}\mathrm{e}^{-z_k\dd_{z_{k-1}}}
\Psi(z_{k-1})\right|_{t_k=1} = \left.\RR(t_k\dd_{t_k}) \Psi(t_k
z_{k-1,k}+z_k)\right|_{t_k=1}
$$
We have
$$
\left[\mathbf{Q}_{-}(u)\Psi\right](z_1\cdots z_N) =
\left.\RR(z_{01}\dd_0)\RR(z_{12}\dd_1)\RR(z_{23}\dd_2)\cdots
\RR(z_{N-1,N}\dd_{N-1})\Psi(z_0,z_1\cdots z_{N-1})\right|_{z_0=z_N}=
$$
$$
= \RR(t_1\dd_{t_1})\RR(t_2\dd_{t_2})\cdots \RR(t_N\dd_{t_N})\cdot
\left.\Psi\left(t_1z_{01}+z_1 , t_2z_{12}+z_2 \cdots
t_Nz_{N-1,N}+z_N\right)\right|_{t_k=1 , z_0=z_n}
$$
which is just the second formula~(\ref{QL2}).

The third formula is the simple consequence of the second one.
We use the integral representation for the operator 
$$
\RR(t_k\dd_{t_k})\left.\Phi(t_k)\right|_{t_k=1} = 
\frac{\Gamma(2\ell)}{\Gamma(\ell+u)\Gamma(\ell-u)}\cdot
\int^{1}_{0}\mathrm{d}\alpha_k (1-\alpha_k)^{\ell-u-1}
\alpha_k^{\ell+u-1}\cdot\alpha_k^{t_k\dd_{t_k}}
\left.\Phi(t_k)\right|_{t_k=1}= 
$$
$$
=\frac{\Gamma(2\ell)}{\Gamma(\ell+u)\Gamma(\ell-u)}\cdot
\int^{1}_{0}\mathrm{d}\alpha_k (1-\alpha_k)^{\ell-u-1}
\alpha_k^{\ell+u-1}\Phi(\alpha_k).
$$
This allows now to obtain~(\ref{QL3}) from~(\ref{QL2}).

The most useful for the proof of the last property 
is the formula~(\ref{QL2}). Let us consider the
action of $\mathbf{Q}_{-}(u)$ on the monomial $z_1^{m_1}\cdots z_N^{m_N}$
$$
\mathbf{Q}_{-}(u) z_1^{m_1}\cdots z_N^{m_N} = \RR(t_1\dd_{t_1})\cdots
\RR(t_N\dd_{t_N})\cdot \left.\left(t_1 z_{N1}+z_1\right)^{m_1}\left( t_2
z_{12}+z_2\right)^{m_2} \cdots \left( t_N
z_{N-1,N}+z_N\right)^{m_N}\right|_{t_k=1}
$$
The left hand side is the sum of monomials $t_1^{k_1}\cdots t_N^{k_N}$ with
polynomials coefficients from $\C\left[z_1\cdots z_N\right]$ so that it is
sufficient to prove that
$$
\RR(t_1\dd_{t_1})\RR(t_2\dd_{t_2})\cdots \RR(t_N\dd_{t_N})\cdot \left.
t_1^{k_1} t_2^{k_2} \cdots t_N^{k_N}\right|_{t_k=1} = \RR(k_1)\RR(k_2)\cdots
\RR(k_n)
$$
is polynomial in $u$. We have
$$
\RR(k)= \frac{\Gamma(2\ell)}{\Gamma(u+\ell)}
\frac{\Gamma\left(k+u+\ell\right)} {\Gamma\left( k+2\ell\right)} = 
\frac{\Gamma(2\ell)}{\Gamma(k+2\ell)}\cdot
(u+\ell)(u+\ell+1)\cdots(u+\ell+k-1)
$$
In the generic situation all is well defined and one obtains a 
polynomial in $u$.
Note that the operator $\mathbf{Q}_{-}(u)$ coincides 
with Q-operator constructed in~\cite{Der1} by using 
Pasquier-Gaudin method~\cite{PG}.
There exists another equivalent representation for the operator 
$\mathbf{Q}_{-}(u)$~\cite{DKM}.  

\section{Back to the inhomogeneous chain:
factorization of the two-parametric operator 
$\mathbf{Q}(u_1|u_2)$}

\setcounter{equation}{0}

There exists the natural generalization of the operators 
$\mathbf{Q}_{\pm}(u)$ to the case of the generic 
inhomogeneous periodic XXX spin chain.
The natural local building blocks for the operators
$\mathbf{Q}_{\pm}(u)$ are the operators $\RR^{\pm}_{k0}$.
\begin{equation}
\mathbf{Q}_{-}(u) = \tr_{V_0} \P_{10}\RR^{-}_{10}(u+\delta_1)\cdots
\P_{N0}\RR^{-}_{N0}(u+\delta_N) \label{Q-}
\end{equation}
\begin{equation}
\mathbf{Q}_{+}(u) = \tr_{V_0} \P_{10}\RR^{+}_{10}(u+\delta_1)\cdots
\P_{N0}\RR^{+}_{N0}(u+\delta_N) \label{Q+}
\end{equation}
where we used the notations
\begin{equation}
\RR^{-}_{k0}(u) = 
\frac{\Gamma\left(2\ell_k\right)}
{\Gamma\left(u+\ell_k\right)}
\cdot
\frac{\Gamma\left(z_{k0}\dd_k+u+\ell_k\right)}
{\Gamma\left(z_{k0}\dd_k+2\ell_k\right)}
\end{equation}
\begin{equation}
\RR^{+}_{k0}(u) = 
\frac{\Gamma\left(1+\ell_k-u\right)}
{\Gamma\left(2\ell_k\right)}
\cdot
\frac{\Gamma\left(z_{0k}\dd_0+2\ell_k\right)}
{\Gamma\left(z_{0k}\dd_0+1+\ell_k-u\right)}.
\end{equation}
Note that explicit expressions for the operator 
$\mathbf{Q}_{-}(u)$ can be obtained by repeating 
step by step all calculations from the previous Section.
\begin{prop}
The operator  $\mathbf{Q}_{-}(u)$ obeys the Baxter's equation
\begin{equation}
\mathbf{Q}_{-}(u)\cdot\mathbf{t}(u) = 
\Delta_{+}(u)\mathbf{Q}_{-}(u+1)+
\Delta_{-}(u)\mathbf{Q}_{-}(u-1)
\ ;\ \Delta_{\pm}(u)=(u+\delta_1\pm\ell_1)\cdots
(u+\delta_N\pm\ell_N)
\label{BQ-}
\end{equation}
and the operator  $\mathbf{Q}_{+}(u)$ obeys 
the Baxter's equation
\begin{equation}
\mathbf{t}(u)\cdot\mathbf{Q}_{+}(u) =
\frac{\Delta_{+}(u-1)\Delta_{-}(u)}{\Delta_{-}(u-1)}
\mathbf{Q}_{+}(u-1) + \Delta_{-}(u)\mathbf{Q}_{+}(u+1)
\label{BQ+}
\end{equation}
\end{prop}
We prove again the first relation only and 
the proof of the second ones is similar.
It is the direct consequence of the triangularity 
relation~(\ref{trianglR-}) and cyclicity of the trace.
Let us choose the first space in~(\ref{trianglR-})
as k-th quantum space and the second 
space as the auxiliary space. Then we have
$$
\P_{k0}\RR^{-}_{k0}(u+\delta_k)\cdot 
\mathrm{L}_k\left(u_k^{+},u_k^{-}\right) =
\mathbf{M}_0\cdot
\left(\begin{array}{cc} u_k^{+}\cdot
\P_{k0}\RR^{-}_{k0}(u+1+\delta_k) & -\P_{k0}\RR^{-}_{k0}(u+\delta_k)\dd_k\\
0 &  u_k^{-}\cdot\P_{k0}\RR^{-}_{k0}(u-1+\delta_k)
\end{array}\right)\ \mathbf{M}^{-1}_0
$$
The triangularity relation for the 
$\RR^{-}$-operator allows to transform the
following product to the triangular form
$$
\P_{10}\RR^{-}_{10}(u+\delta_1)\cdots
\P_{N0}\RR^{-}_{N0}(u+\delta_N)\cdot
\mathrm{L}_1\left(u_1^{+},u_1^{-}\right)
\cdots
\mathrm{L}_N\left(u_N^{+},u_N^{-}\right) =
$$
$$
=\P_{10}\RR^{-}_{10}(u+\delta_1)\mathrm{L}_1\left(u_1^{+},u_1^{-}\right)
\cdots
\P_{N0}\RR^{-}_{N0}(u+\delta_N)\mathrm{L}_N\left(u_N^{+},u_N^{-}\right)=
$$
$$
= \mathbf{M}_0
\left(\begin{array}{cc} u^{+}_1\P_{10}\RR^{-}_{10}(u+1+\delta_1)& *** \\
0 & u^{-}_1\P_{10}\RR^{-}_{10}(u-1+\delta_1)
\end{array}\right)
\cdots
$$
$$
\cdots\left(\begin{array}{cc} u^{+}_N
\P_{N0}\RR^{-}_{N0}(u+1+\delta_N)& *** \\
0 & u^{-}_N\P_{N0}\RR^{-}_{N0}(u-1+\delta_N)
\end{array}\right)
\mathbf{M}^{-1}_0 
$$
To derive the Baxter's equation it remains to calculate the traces in the
auxiliary spaces $V_0$ and $\C^2$ using the following simple statement.
Let operators $\mathbb{A}_{ij}$ act in the tensor product $V_0\otimes
V_1\cdots\otimes V_n$. Then we have
$$
\tr_{V_0}\tr \left(\begin{array}{cc}
1 & 0 \\
z_0 & 1\end{array}\right) \left(\begin{array}{cc}
\mathbb{A}_{11} & \mathbb{A}_{12} \\
0 & \mathbb{A}_{22}\end{array}\right) \left(\begin{array}{cc}
1 & 0 \\
-z_0 & 1\end{array}\right) =
\tr_{V_0}\mathbb{A}_{11}+\tr_{V_0}\mathbb{A}_{22}.
$$
The proof is straightforward and uses the cyclic property of the trace only
$$
 \tr_{V_0}\mathbb{A}_{12}\cdot z_0 =
 \tr_{V_0} z_0\cdot \mathbb{A}_{12}.
$$
Next we consider the commutation relations between the operators 
$\mathbf{Q}_{\pm}(u)$ and the transfer matrix $\mathbf{t}(u)$.
\begin{prop}
The operators $\mathbf{Q}_{\pm}(u)$ have the following 
commutation relations with the transfer matrix~(\ref{t})
$$
\mathbf{Q}_{-}\left(\lambda\right)\cdot 
\tr \mathrm{L}_{1}\left(u_1^{+},u_1^{-}\right)
\mathrm{L}_{2}\left(u_2^{+},u_2^{-}\right)\cdots
\mathrm{L}_{n}\left(u_n^{+},u_n^{-}\right) =
$$
\begin{equation}
= \tr \mathrm{L}_{1}\left(u_2^{+},u_1^{-}\right)
\mathrm{L}_{2}\left(u_3^{+},u_2^{-}\right)\cdots
\mathrm{L}_{n}\left(u_1^{+},u_n^{-}\right)\cdot 
\mathbf{Q}_{-}\left(\lambda\right),
\label{Q-LL}
\end{equation}
$$
\tr \mathrm{L}_{1}\left(u_1^{+},u_1^{-}\right)
\mathrm{L}_{2}\left(u_2^{+},u_2^{-}\right)\cdots
\mathrm{L}_{n}\left(u_n^{+},u_n^{-}\right)\cdot 
\mathbf{Q}_{+}\left(\lambda\right) =
$$
\begin{equation}
= \mathbf{Q}_{+}\left(\lambda\right)\cdot \tr
\mathrm{L}_{1}\left(u_1^{+},u_n^{-}\right)
\mathrm{L}_{2}\left(u_2^{+},u_1^{-}\right)\cdots
\mathrm{L}_{n}\left(u_n^{+},u_{n-1}^{-}\right).
\label{Q+LL}
\end{equation}
\end{prop}
Let us prove the equation~(\ref{Q-LL}) for
example. We start from the commutation relation 
$$
\P_{10}\RR^{+}_{10}(u+\delta_1-\lambda_{-})\cdots
\P_{N0}\RR^{+}_{N0}(u+\delta_N-\lambda_{-})\cdot
\mathrm{L}_{1}\left(u_1^{+},u_1^{-}\right)
\mathrm{L}_{2}\left(u_2^{+},u_2^{-}\right)\cdots
\mathrm{L}_{n}\left(u_n^{+},u_n^{-}\right)\cdot
\mathrm{L}_{0}\left(\lambda_{+},\lambda_{-}\right) =
$$
$$
= \mathrm{L}_{0}\left(u_1^{+},\lambda_{-}\right)\cdot
\mathrm{L}_{1}\left(u_2^{+},u_1^{-}\right)
\mathrm{L}_{2}\left(u_3^{+},u_2^{-}\right)\cdots
\mathrm{L}_{n}\left(\lambda_{+},u_n^{-}\right)\cdot
\P_{10}\RR^{+}_{10}(u+\delta_1-\lambda_{-})\cdots
\P_{N0}\RR^{+}_{N0}(u+\delta_N-\lambda_{-})
$$
which is derived directly from the defining relation~(\ref{F2}) 
for the operator $\RR^{-}$.
Next we put
$\lambda_{+} = u_1^{+}$ and multiply both sides of the equation by the
$\mathrm{L}^{-1}_{0}\left(u_1^{+};\lambda_{-}\right)$. It remains to
calculate the traces in the auxiliary spaces $V_0$ and $\C^2$ using the
following simple statement.
Let operators $\mathbb{A}_{ij}$ act in the tensor product $V_0\otimes
V_1\cdots\otimes V_n$.
\begin{equation}
\tr_{V_0}\tr \mathrm{L}_0\left(u^{+};u^{-}\right)
\cdot\left(\begin{array}{cc}
\mathbb{A}_{11} & \mathbb{A}_{12} \\
\mathbb{A}_{21} & \mathbb{A}_{22}\end{array}\right)\cdot
\mathrm{L}^{-1}_0\left(u^{+};u^{-}\right) =
\tr_{V_0}\mathbb{A}_{11}+\tr_{V_0}\mathbb{A}_{22}.
\end{equation}
The proof is straightforward and uses the cyclic property of the trace.
After all one obtains the commutation relation~(\ref{Q-LL}) for $\lambda
= u-\lambda_{-}$ which is equivalent~(\ref{Q-LL}) due to arbitrariness of
$\lambda_{-}$.

It is possible to build from the operators $\mathbf{Q}_{\pm}$ 
some composite operator commuting with 
transfer matrix $\mathbf{t}(u)$.
The operator
$
\mathbf{Q}_{+}\left(\lambda\right)\cdot\P
\cdot\mathbf{Q}_{-}\left(\mu\right),
$
where $\P$ is the operator of cyclic shift
$$
\P \psi\left(z_1,z_2\cdots z_n\right) = \psi\left(z_2,z_3\cdots z_n,
z_1\right)
$$
commutes with the generic transfer matrix~(\ref{t})
$$
\mathbf{Q}_{+}\left(\lambda\right)\cdot\P
\cdot\mathbf{Q}_{-}\left(\mu\right)\cdot \tr
\mathrm{L}_{1}\left(u_1^{+},u_1^{-}\right)
\mathrm{L}_{2}\left(u_2^{+},u_2^{-}\right)\cdots
\mathrm{L}_{n}\left(u_n^{+},u_n^{-}\right) = 
$$
$$
= \tr \mathrm{L}_{1}\left(u_1^{+},u_1^{-}\right)
\mathrm{L}_{2}\left(u_2^{+},u_2^{-}\right)\cdots
\mathrm{L}_{n}\left(u_n^{+},u_n^{-}\right)\cdot 
\mathbf{Q}_{+}\left(\lambda\right)\cdot\P
\cdot\mathbf{Q}_{-}\left(\mu\right)
$$
The commutativity follows from the commutation
relations~(\ref{Q-LL}) and~(\ref{Q+LL})
$$
\mathbf{Q}_{+}\left(\lambda\right)\cdot\P\cdot
\mathbf{Q}_{-}\left(\mu\right)\cdot \tr
\mathrm{L}_{1}\left(u_1^{+},u_1^{-}\right)
\mathrm{L}_{2}\left(u_2^{+},u_2^{-}\right)\cdots
\mathrm{L}_{n}\left(u_n^{+},u_n^{-}\right) =
$$
$$
= \mathbf{Q}_{+}\left(\lambda\right)\cdot\P\cdot\tr
\mathrm{L}_{1}\left(u_2^{+},u_1^{-}\right)
\mathrm{L}_{2}\left(u_3^{+},u_2^{-}\right)\cdots
\mathrm{L}_{n}\left(u_1^{+},u_n^{-}\right)\cdot 
\mathbf{Q}_{-}\left(\mu\right)=
$$
$$
= \mathbf{Q}_{+}\left(\lambda\right)\cdot\tr
\mathrm{L}_{2}\left(u_2^{+},u_1^{-}\right)
\mathrm{L}_{3}\left(u_3^{+},u_2^{-}\right)\cdots
\mathrm{L}_{1}\left(u_1^{+},u_n^{-}\right)\cdot 
\P \cdot\mathbf{Q}_{-}\left(\mu\right)=
$$
$$
= \mathbf{Q}_{+}\left(\lambda\right)\cdot\tr
\mathrm{L}_{1}\left(u_1^{+},u_n^{-}\right)
\mathrm{L}_{2}\left(u_2^{+},u_1^{-}\right)\cdots
\mathrm{L}_{n}\left(u_n^{+},u_{n-1}^{-}\right)\cdot
\P\cdot\mathbf{Q}_{-}\left(\mu\right)=
$$
$$
=\tr \mathrm{L}_{1}\left(u_1^{+},u_1^{-}\right)
\mathrm{L}_{2}\left(u_2^{+},u_2^{-}\right)\cdots
\mathrm{L}_{n}\left(u_n^{+},u_n^{-}\right)\cdot
\mathbf{Q}_{+}\left(\lambda\right)\cdot\P\cdot 
\mathbf{Q}_{-}\left(\mu\right)
$$
The operator $\mathbf{Q}_{+}\left(\lambda\right)\cdot\P\cdot 
\mathbf{Q}_{-}\left(\mu\right)$ obeys both Baxter equations 
with respect to each parameter. 
It is a simple consequence of Baxter's equations for the operators 
$\mathbf{Q}_{+}(\lambda)$ and $\mathbf{Q}_{-}(\mu)$.

After all it is evident that in the case 
of homogeneous spin chain: $\delta_k=0$ and $\ell_k=\ell$ 
we have the points of degeneracy $\lambda = 1-\ell$ 
and $\mu =\ell$ where operators $\mathbf{Q}_{\pm}$ 
reduced to the operator $\P$:
$\mathbf{Q}_{+}\left(1-\ell\right)= \P$ and 
$\mathbf{Q}_{-}\left(\ell\right)= \P$.
So that we have 
$$
\mathbf{Q}_{+}\left(1-\ell\right)\cdot\P\cdot 
\mathbf{Q}_{-}\left(\mu\right)=\mathbf{Q}_{-}\left(\mu\right)  
\ ;\ \mathbf{Q}_{+}\left(\lambda\right)\cdot\P\cdot 
\mathbf{Q}_{-}\left(\ell\right) = 
\mathbf{Q}_{+}\left(\lambda\right)
$$
We see that the composite operator 
$\mathbf{Q}_{+}\left(u_1\right)\cdot\P\cdot 
\mathbf{Q}_{-}\left(u_2\right)$ has the same properties 
as the two-parametric operator $\mathbf{Q}\left(u_1|u_2\right)$.
Of course this is not accidental because these 
operators indeed coincide.
\begin{prop}
\begin{equation}
\mathbf{Q}\left(u_1|u_2\right)=
\mathbf{Q}_{+}\left(u_1\right)\cdot\P\cdot 
\mathbf{Q}_{-}\left(u_2\right)
\label{factor+-}
\end{equation}
\end{prop}
The direct proof which we have is rather technical and 
does not illuminate the origin of this factorization.
For brevity we omit the proof and hope that this 
factorization looks sufficiently natural.

Finally we consider the commutation relations 
between operators $\mathbf{Q}_{+}\left(u_1\right)$ and 
$\mathbf{Q}_{-}\left(u_2\right)$.
\begin{prop}
The operators 
$\mathbf{Q}_{+}\left(u_1\right)$ and 
$\mathbf{Q}_{+}\left(u_1\right)$ have the following 
commutation relations
\begin{equation}
\mathbf{Q}_{+}\left(u_1\right)
\cdot\P\cdot\mathbf{Q}_{-}\left(u_2\right)
\cdot\mathbf{Q}_{+}\left(v_1\right)=
\mathbf{Q}_{+}\left(v_1\right)
\cdot\P\cdot\mathbf{Q}_{-}\left(u_2\right)
\cdot\mathbf{Q}_{+}\left(u_1\right)
\label{Q+}
\end{equation}
\begin{equation}
\mathbf{Q}_{-}\left(u_2\right)
\cdot\mathbf{Q}_{+}\left(v_1\right)\cdot\P\cdot
\mathbf{Q}_{-}\left(v_2\right)= 
\mathbf{Q}_{-}\left(v_2\right)
\cdot\mathbf{Q}_{+}\left(v_1\right)\cdot\P\cdot
\mathbf{Q}_{-}\left(u_2\right)
\label{Q-}
\end{equation}
\end{prop}
These commutation relations are the consequence of the 
commutation relations for the two-parametric operators 
\begin{equation}
\mathbf{Q}(u_1|u_2)\cdot\mathbf{Q}(v_1|v_2)=
\mathbf{Q}(v_1|u_2)\cdot\mathbf{Q}(u_1|v_2)
\label{1}
\end{equation}
\begin{equation}
\mathbf{Q}(u_1|u_2)\cdot\mathbf{Q}(v_1|v_2)=
\mathbf{Q}(u_1|v_2)\cdot\mathbf{Q}(v_1|u_2)
\label{2}
\end{equation}
and the factorization of the two-parametric operator
$$
\mathbf{Q}\left(u_1|u_2\right)=
\mathbf{Q}_{+}\left(u_1\right)\cdot\P\cdot 
\mathbf{Q}_{-}\left(u_2\right).
$$
We shall prove the commutation relations~(\ref{1}) and~(\ref{2}).  
These commutation relations follow from some local 
three-term relations which are similar 
to the Yang-Baxter relation for the $\R$-operators.  
To derive the needed relation we shall proceed in  
close analogy with the well known derivation of the Yang-Baxter 
relation from the defining equation
$$
\check{\R}_{12}(u_{+},u_-|v_{+},v_-)
\mathrm{L}_{1}(u_+,u_-)\mathrm{L}_{2}(v_+,v_-)=
\mathrm{L}_{1}(v_+,v_-)\mathrm{L}_{2}(u_+,u_-)
\check{\R}_{12}(u_{+},u_-|v_{+},v_-)
$$
We recall that the commutativity of the diagram

\vspace{5mm} \unitlength 0.8mm \linethickness{0.4pt}
\begin{picture}(152.67,83.33)
\put(30.00,50.00){\makebox(0,0)[cc]
{$\mathrm{L}_1(u_{+},u_-)\mathrm{L}_2(v_{+},v_-)L_3(w_{+},w_-)$}}
\put(32.00,55.00){\vector(0,1){20.00}}
\put(57.00,65.00){\makebox(0,0)[cc]
{$\check{\R}_{12}(u_{+},u_-|v_{+},v_-)$}}
\put(30.00,80.00){\makebox(0,0)[cc]
{$\mathrm{L}_1(v_{+},v_-)\mathrm{L}_2(u_{+},u_-)L_3(w_{+},w_-)$}}
\put(65.00,80.00){\vector(1,0){40.00}}
\put(85.00,85.00){\makebox(0,0)[cc]
{$\check{\R}_{23}(u_{+},u_-|w_+,w_-)$}}
\put(140.00,80.00){\makebox(0,0)[cc]
{$\mathrm{L}_1(v_{+},v_-)\mathrm{L}_2(w_{+},w_-)L_3(u_{+},u_-)$}}
\put(140.00,75.00){\vector(0,-1){20.00}}
\put(117.00,65.00){\makebox(0,0)[cc]
{$\check{\R}_{12}(v_+,v_-|w_{+},w_-)$}}
\put(140.00,50.00){\makebox(0,0)[cc]
{$\mathrm{L}_1(w_{+},w_-)\mathrm{L}_2(v_{+},v_-)L_3(u_{+},u_-)$}}
\put(32.00,45.00){\vector(0,-1){20.00}}
\put(57.00,35.00){\makebox(0,0)[cc]
{$\check{\R}_{23}(v_{+},v_{-}|w_+,w_-)$}}
\put(30.00,20.00){\makebox(0,0)[cc]
{$\mathrm{L}_1(u_{+},u_-)\mathrm{L}_2(w_{+},w_-)L_3(v_{+},v_-)$}}
\put(65.00,20.00){\vector(1,0){40.00}}
\put(85.00,15.00){\makebox(0,0)[cc]
{$\check{\R}_{12}(u_+,u_{-}|w_+,w_{-})$}}
\put(140.00,20.00){\makebox(0,0)[cc]
{$\mathrm{L}_1(w_{+},w_-)\mathrm{L}_2(u_{+},u_-)L_3(v_{+},v_-)$}}
\put(140.00,25.00){\vector(0,1){20.00}}
\put(117.00,35.00){\makebox(0,0)[cc]
{$\check{\R}_{23}(u_{+},u_-|v_{+},v_-)$}}
\end{picture}

results in the Yang-Baxter equation for $\R$-operators
$$
\check{\R}_{12}(v_{+},v_-|w_{+},w_-)\check{\R}_{23}(u_{+},u_-|w_+,w_-)
\check{\R}_{12}(u_+,u_-|v_{+},v_-) =
$$
$$
=\check{\R}_{23}(u_{+},u_{-}|v_+,v_-)\check{\R}_{12}(u_+,u_{-}|w_+,w_{-})
\check{\R}_{23}(v_{+},v_-|w_{+},w_-).
$$
In the same way the commutativity of the diagram

\vspace{5mm} \unitlength 0.8mm \linethickness{0.4pt}
\begin{picture}(152.67,83.33)
\put(30.00,50.00){\makebox(0,0)[cc]
{$\mathrm{L}_1(u_{+},u_-)\mathrm{L}_2(v_{+},v_-)L_3(w_{+},w_-)$}}
\put(32.00,55.00){\vector(0,1){20.00}}
\put(57.00,65.00){\makebox(0,0)[cc]
{$\check{\R}_{12}(u_{+},u_-|v_+,v_{-})$}}
\put(30.00,80.00){\makebox(0,0)[cc]
{$\mathrm{L}_1(v_{+},v_-)\mathrm{L}_2(u_{+},u_-)L_3(w_{+},w_-)$}}
\put(65.00,80.00){\vector(1,0){40.00}}
\put(85.00,85.00){\makebox(0,0)[cc]
{$\check{\R}_{23}(u_{+},u_-|w_+,w_-)$}}
\put(140.00,80.00){\makebox(0,0)[cc]
{$\mathrm{L}_1(v_{+},v_-)\mathrm{L}_2(w_{+},w_-)L_3(u_{+},u_-)$}}
\put(140.00,75.00){\vector(0,-1){20.00}}
\put(117.00,65.00){\makebox(0,0)[cc]
{$\RR^{-}_{12}(v_+,v_-|w_-)$}}
\put(140.00,50.00){\makebox(0,0)[cc]
{$\mathrm{L}_1(v_{+},w_-)\mathrm{L}_2(w_{+},v_-)L_3(u_{+},u_-)$}}
\put(32.00,45.00){\vector(0,-1){20.00}}
\put(57.00,35.00){\makebox(0,0)[cc]
{$\RR^{-}_{23}(v_{+},v_{-}|w_-)$}}
\put(30.00,20.00){\makebox(0,0)[cc]
{$\mathrm{L}_1(u_{+},u_-)\mathrm{L}_2(v_{+},w_-)L_3(w_{+},v_-)$}}
\put(65.00,20.00){\vector(1,0){40.00}}
\put(85.00,15.00){\makebox(0,0)[cc]
{$\check{\R}_{12}(u_+,u_{-}|v_+,w_{-})$}}
\put(140.00,20.00){\makebox(0,0)[cc]
{$\mathrm{L}_1(v_{+},w_-)\mathrm{L}_2(u_{+},u_-)L_3(w_{+},v_-)$}}
\put(140.00,25.00){\vector(0,1){20.00}}
\put(117.00,35.00){\makebox(0,0)[cc]
{$\check{\R}_{23}(u_{+},u_-|w_+,v_{-})$}}
\end{picture}

results in three term relation for one $\RR^{-}$- and two
$\check{\R}$-operators
\begin{equation}
\RR^{-}_{12}(v_+,v_-|w_-)\check{\R}_{23}(u_{+},u_-|w_+,w_-)
\check{\R}_{12}(u_{+},u_-|v_+,v_{-})
=
\label{R-RR}
\end{equation}
$$
=\check{\R}_{23}(u_{+},u_-|w_+,v_{-})
\check{\R}_{12}(u_+,u_{-}|v_+,w_{-})
\RR^{-}_{23}(v_{+},v_{-}|w_-)
$$
In a similar way one obtains the relation for one
$\RR^{+}$- and two $\check{\R}$-operators
\begin{equation}
\RR^{+}_{12}(v_+|w_+,w_-)\check{\R}_{23}(u_{+},u_-|w_+,w_-)
\check{\R}_{12}(u_{+},u_-|v_+,v_{-}) =
\label{R+RR}
\end{equation}
$$
=\check{\R}_{23}(u_{+},u_-|v_+,w_{-})\check{\R}_{12}(u_+,u_{-}|w_+,v_{-})
\RR^{+}_{23}(v_{+}|w_{+},w_-).
$$
First of all after multiplication with the permutation operators 
the three term relation~(\ref{R-RR}) can be rewritten as follows
$$
\P_{23}\RR^{-}_{23}(v_+,v_-|w_-)\R_{13}(u_{+},u_-|w_+,w_-)
\R_{12}(u_{+},u_-|v_+,v_{-}) =
$$
$$
=\R_{12}(u_{+},u_-|w_+,v_{-})\R_{13}(u_+,u_{-}|v_+,w_{-})
\P_{23}\RR^{-}_{23}(v_{+},v_{-}|w_-)
$$
Next we choose the first space $\mathrm{V}_k$, 
the second space $\mathrm{V}_0$ and the 
third space $\mathrm{V}_{0^{\prime}}$.
We have
$$
\P_{00^{\prime}}\RR^{-}_{00^{\prime}}(v_+,v_-|w_-)
\R_{k0^{\prime}}(u_{+},u_-|w_+,w_-)
\R_{k0}(u_{+},u_-|v_+,v_{-}) =
$$
$$
=\R_{k0}(u_{+},u_-|w_+,v_{-})\R_{k0^{\prime}}(u_+,u_{-}|v_+,w_{-})
\P_{00^{\prime}}\RR^{-}_{00^{\prime}}(v_{+},v_{-}|w_-)
$$
These local relations imply in the standard way the 
following commutation relation 
$$
\mathbf{Q}(u_{1}|u_2)\mathbf{Q}(v_1|v_2) =
\mathbf{Q}(v_{1}|u_2)\mathbf{Q}(u_1|v_2)
$$
In a similar way choosing the first space $\mathrm{V}_k$, 
the second space $\mathrm{V}_0$ and the 
third space $\mathrm{V}_{0^{\prime}}$ we rewrite the equation
$$
\RR^{+}_{12}(v_+|w_+,w_-)\check{\R}_{23}(u_{+},u_-|w_+,w_-)
\check{\R}_{12}(u_{+},u_-|v_+,v_{-}) =
$$
$$
=\check{\R}_{23}(u_{+},u_-|v_+,w_{-})
\check{\R}_{12}(u_+,u_{-}|w_+,v_{-})
\RR^{+}_{23}(v_{+}|w_{+},w_-)
$$
in the form
$$
\P_{00^{\prime}}\RR^{+}_{00^{\prime}}(v_+|w_+,w_-)
\R_{k0^{\prime}}(u_{+},u_-|w_+,w_-)
\R_{k0}(u_{+},u_-|v_+,v_{-}) =
$$
$$
=\R_{k0}(u_{+},u_-|v_+,v_{-})\R_{k0^{\prime}}(u_{+},u_-|w_+,w_-)
\P_{00^{\prime}}\RR^{+}_{00^{\prime}}(v_+|w_+,w_-).
$$
These local relations imply in the standard way the 
following commutation relation  
$$
\mathbf{Q}(u_{1}|u_2)\mathbf{Q}(v_1|v_2) =
\mathbf{Q}(u_{1}|v_2)\mathbf{Q}(v_1|u_2).
$$

\section{Conclusions}

Using the unversal R-matrices $\R_{k0}(u)$ as 
building blocks it is possible to construct 
the two-parametric Baxter's Q-operator in the case of 
generic inhomogeneous periodic XXX spin chain 
$$
\mathbf{Q}(u_1|u_2) = 
\tr_{\mathrm{V}_{0}}\R_{10}(u+\delta_1)\R_{20}(u+\delta_2)\cdots 
\R_{N0}(u+\delta_N).
$$
This operator is factorized on the product 
of simpler operators $\mathbf{Q}_{+}(u_1)$ 
and $\mathbf{Q}_{-}(u_2)$
$$
\mathbf{Q}(u_1|u_2) = \mathbf{Q}_{+}(u_1)
\cdot\P\cdot\mathbf{Q}_{-}(u_2).  
$$
In  the general case of inhomogeneous spin 
chain the operators $\mathbf{Q}_{+}(u_1)$ and 
$\mathbf{Q}_{-}(u_2)$ are not $s\ell(2)$-invariant.
There are nontrivial commutation relations between 
$\mathbf{Q}_{+}$ and $\mathbf{Q}_{-}$ and 
between $\mathbf{Q}_{\pm}$ and transfer 
matrix $\mathrm{t}(u)$.
In the special case of homogeneous spin chain 
the operators $\mathbf{Q}_{+}(u)$,$\mathbf{Q}_{-}(u)$ and 
$\mathbf{t}(u)$ become commuting operators and  
the $s\ell(2)$-invariance of the operators 
$\mathbf{Q}_{+}(u)$ and $\mathbf{Q}_{-}(u)$ is restored. 
The operators $\mathbf{Q}_{+}(u)$ and $\mathbf{Q}_{-}(u)$  
have all properties of Baxter's Q-operator in the special case 
of homogeneous spin chain but in the case of 
inhomogeneous spin only the composite operator $\mathbf{Q}_{+}(u_1)
\cdot\P\cdot\mathbf{Q}_{-}(u_2)$ has the needed properties.

The factorization of the general 
R-matrix can be generalized to the more complicated 
situation when the symmetry algebra is $s\ell(3)$~\cite{Der}. 
We hope that there exists the similar construction 
of the Baxter's Q-operator in this case.

\section{Acknowledgments}

I would like to thank R.Kirschner, G.Korchemsky, P.Kulish,
A.Manashov, E.Sklyanin and~V.Tarasov for the stimulating
discussions and critical remarks on the different stages of
this work. I thank also the Center of theoretical 
science of Leipzig University for hospitality in the 
stage of writing up this paper. 
This work was supported by the grant 03-01-00837
of the Russian Foundation for Fundamental Research and 
partially by a travel grant of Deutsche Forschungsgemeinschaft.

\end{document}